

\magnification=1200
\input amstex
\documentstyle{amsppt}
\define\CDdashright#1#2{&\,\mathop{\longdashrightarrow}\limits^{#1}_{#2}\,&}
\define\CDdashleft#1#2{&\,\mathop{\longdashleftarrow}\limits^{#1}_{#2}\,&}

\def\Z{\Bbb Z}
\def\P{\Bbb P}

\def\C{\Bbb C}

\def\Til#1{\widetilde{#1}}

\def\Pic{\text{Pic}}

\def\tablerule{\noalign{\hrule}}
\def\siltable#1.{$$\text{
\vbox{\tabskip=0pt \offinterlineskip
\halign to300pt{\strut##& ##\tabskip=1em plus2em&
  \hfil##\hfil& \vrule##&
  \hfil##\hfil& \vrule##&
  \hfil##\hfil& \vrule##&
  \hfil##\hfil& \vrule##&
  \hfil##\hfil& ##\tabskip=0pt\cr
#1}}}$$}

\pagewidth{12.45 cm}\pageheight{20.35 cm}
\topmatter
\title A blow-up construction and graph coloring\endtitle
\author Paolo Aluffi \endauthor
\affil Mathematics Department, Florida State University\endaffil
\address Tallahassee, FL 32306
\endaddress
\endtopmatter
\document


\heading \S 0. Introduction\endheading

In this note we construct a non@-singular algebraic variety $V_G$
encoding the incidence information of a simple graph $G$, by a
sequence of blow@-ups of a projective space along suitable linear
subspaces. The aim is to translate into the geometry of $V_G$ the
combinatorial information about $G$; we find that this can be done
with surprising ease and efficiency.

For example, we prove that the {\it chromatic polynomial\/} of the
graph---that is, the polynomial giving for each $m>0$ the number of
ways in which $G$ can be colored using $m$ colors, so that no two
adjacent vertices are assigned the same color---is (up to a power of
the variable) the intersection product of a fixed class $\gamma$ in
$A_1 V_G$ with a polynomial $S(t)$ in $\Pic V_G[t]$: the class is
defined as the Poincar\'e dual of the pull@-back of the hyperplane
class, with respect to a natural basis of $\Pic V_G$, and $S(t)$ is
also easily defined as a combination of the exceptional divisors
arising in the blow@-up construction. In \S 1 we describe the
construction for graphs and state the above result precisely (but with
no proofs), as a sales pitch for the rest of the paper, which examines
the construction more carefully and gives deeper---but necessarily
more technical---results.

In fact the right level of generality to perform our construction is
that of `combinatorial geometries which are projectively
coordinatizable over some field'; for short (and a little improperly)
we will refer to these as {\it matroids\/}. Our construction can be
performed starting from any (loopless) matroid embedded in a
projective space, and specializes to the one in \S1 for the {\it cycle
matroid\/} of a graph. We give this more general construction in \S 2:
roughly, the variety of a matroid is obtained by blowing@-up the
ambient projective space along the flats of the matroid, in order of
increasing dimension. We prove the above result concerning the
chromatic polynomial of a graph by showing that the {\it
characteristic polynomial\/} of a matroid equals the intersection
product of a fixed 1@-class by a suitable polynomial $S(t)$ in the
$\Pic$ of its variety. A question that then arises naturally regards
the positivity of $S(m)$ for a given $m$ and a given class of
matroids: we determine a large class (including cycle matroids of
graphs) for which a close relative $\overline S(m)$ of $S(m)$ is
generated by global sections for all positive $m$.

To support the point that our construction may offer a new angle on
the theory of characteristic polynomials of matroids, in \S 3 we give
`geometric proofs' of a few basic results on these (our source of
examples here is Zaslavsky's contribution to~\cite{W2}). The
deletion@-contraction rule and Stanley's `modular factorization
theorem' for example follow easily from the functoriality of the
construction. Most likely these proofs could be translated word by
word into standard combinatorial proofs; our point here is that our
arguments are suggested by `algebro@-geometrical intuition', and the
hope is that this could lead to a fresh approach to the combinatorics.
Also, we hope \S3 will help to advertise this beautiful branch of
combinatorics among the geometers.

Our favorite example of the interplay between the two fields is the
following: if we were to hand our construction to a random algebraic
geometer, and asked to provide us with an interesting numerical
invariant of these objects, she would likely propose the intersection
product of the canonical divisor (which is the first place where to
look for an invariant) with the above class $\gamma$ (dual to the
pull@-back of the hyperplane class, thus a priori defined for all
varieties produced by the construction).  The result would essentially
be, as we show in \S 3, Crapo's {\it Beta invariant\/} of the matroid;
the basic properties of this latter (like additivity, or vanishing for
disconnected matroids) all follow from the adjunction formula for the
canonical divisor.

One feature of our construction is that it produces an infinite
tower of varieties, rather than a single one: the construction depends
on a starting $\P^n$ in which the matroid is embedded, and we get a
non@-singular variety $V^n$ of dimension $n$ for each $n$ strictly
larger than the rank of the matroid. In addition, each $V^n$ is
naturally embedded as a divisor in $V^{n+1}$, in a way compatible with
the construction: for example, the divisor $S(m)$ on $V^n$ is the
restriction of the corresponding divisor on $V^{n+1}$, etc. The facts
discussed in the first three sections hold uniformly for each variety
in the sequence, so we may choose one arbitrarily if we wish.  We
think however that interesting information can be extracted from the
whole tower: one such facts is observed in \S 4.  For simplicity,
assume the matroid to be regular (for example, graphical) and consider
the rational maps $V^n \dashrightarrow \P^N$ defined by $S(m)$. Define
$d(m,n)$ to be the degree of the (closure of the) image of this map as
a cycle of dimension $n$. These numbers are invariants of the starting
matroid which, we argue, encode interesting information. $d(m,n)$ is
hard to compute in general (this is almost always the case for the
degree of the image of a rational map!); specific examples can however
be worked out. Here is a table of $d(m,n)$ for a few small values of
$m,n$, for the varieties constructed starting from the complete graph
on three vertices (these entries and the table in \S4 were checked with
Schubert \cite{K-S}):
\siltable
&&$d(m,n)$&& $m=2$ && 3 && 4 && 5 &\cr\tablerule
&& $n=3$  && 42 && 644 && 3888 && 15216 &\cr\tablerule
&&   4    && 210 && 6312 && 64746 && 388704 &\cr\tablerule
&&   5    && 930 && 58312 && 1045476 && 9756192 &\cr\tablerule
&&   6    && 3906 && 529244 && 16764894 && 244093680 &\cr\tablerule
&&   7    && 16002 && 4776396 && 268386264 && 6103281168 &\cr.

And here is the general result given in \S 4:
\proclaim{Theorem} If $n$ is prime and greater than the rank of the
matroid, then
$$d(m,n) \equiv p(m) \pmod n\quad,$$
where $p(m)$ is the characteristic polynomial of the
matroid.\endproclaim

\noindent For example, $p(m)=(m-1)(m-2)$ for the complete graph on
three vertices, and e.g., $6103281168\equiv 4\cdot 3 \pmod 7$.

We note that, by this result, the statement of the celebrated
four@-color@-theorem tranlates into: {\it For a planar graph with $N$
vertices, there exists a prime $n\ge N$ such that $d(4,n)$ is not a
multiple of $n$.}

The above table will immediately convince the reader that it is {\it
not\/} true that $d(m,n)\equiv p(m)\pmod n$ for all $n$.

The numbers $d(m,n)$ above can also be defined without ever leaving
the original projective space from which the construction starts: they
can be written in terms of the Segre classes of specific schemes
supported on a linear subspace of the projective space. A congruence
formula similar to the above can then be written for the
zero@-dimensional term of these Segre classes; see \S4 for a precise
statement.\vskip 12pt

Translating coloring problems in terms of projective geometry is not a
new idea: the `critical problem' (\cite{C-R}, Chapter 16) is the
foremost such construction. We also know of a different and more
algebro-geometric interpretation of these problems due to R. Miranda
(\cite{M}; see also \cite{C-M}). A feature common to the critical
problem and Miranda's approach is that both work by coordinatizing the
relevant combinatorial geometry over a {\it finite\/} field, which in
a sense keeps track of the number of colors used. Our construction has
a different flavor in that it is performed in any characteristic over
which the relevant matroid can be embedded (for example over $\C$;
graphical matroids can be embedded in any characteristic); different
colorings correspond to different divisors within this one
construction. Of course we would be very interested in learning about
relations between our construction and Miranda's or the critical
problem.

Granted, we offer no new coloring theorem here. One missing ingredient
is an algebro-geometric tool to tell a priori when a variety $V_G$ as
above does in fact arise from a planar graph as per our construction:
the next natural step to take in the program is a suitable translation
of Kuratowski's theorem in this language.\vskip 12pt

A note about our references: we draw most of our combinatorics
know@-how from Crapo and Rota's `Combinatorial geometries'
(\cite{C-R}) and from the excellent contributions of Brylawski and
Zaslavsky to \cite{W1}, \cite{W2}. We found these references extremely
helpful for their thoroughness and accessibility to the complete
outsider, of which this writer is a perfect specimen.

Finally, a note for the hasty reader: the reader who feels confortable
with matroids can safely skip \S 1, which simply specializes the
construction to graphs. Also, \S4 can be read independently of
\S3.\vskip 12pt


\heading \S 1. The chromatic polynomial as an intersection
product\endheading

Let $G$ be a graph (loopless and with no parallel edges). Place the
vertices of $G$ at linearly independent points of a projective space
$\P^n$ (over any algebraically closed field), and draw for each edge
the line joining the corresponding vertices. Intersecting the
resulting reducible curve with a general hyperplane gives a
configuration of points $e_k$ (ordered in any fashion), each
corresponding to an edge of the graph, which is the starting point of
our construction: in \S 2 we will study more generally the
construction obtained by starting with any finite collection of points
in a projective space. Our goal is to extract information from the
linear dependence of the points $e_k$; the above is the standard way
to embed in a projective space the `cycle matroid' corresponding to
the graph. The (point corresponding to an) edge $e$ is in the subspace
spanned by edges $e_1,\dots,e_d$ if and only if $e$ joins vertices in
one connected component of the subgraph of $G$ determined by
$e_1,\dots,e_d$ (cf.~e.g.~\cite{W1}, p.19, or \cite{C-R}, chapter 6).
For example, three of the $e_k$'s are collinear in $\P^n$ precisely if
the corresponding edges form a circuit in $G$.

Now for the construction. Consider all dimension@-$d$ subspaces
$x^d_r$ spanned by the $e_k$ in $\P^n$, listed by dimension and
otherwise in any order: so in particular the $x^0_r$'s are simply the
$e_k$'s. Also, consider the subspaces $y^d_r$ obtained by intersecting
collections of the $x$'s, provided these do not appear already in the
list of the $x$'s. Observe that the family of subspaces of $\P^n$ thus
obtained is closed with respect to intersection.

Let $V_0=\P^n$, and inductively let $V_{d+1}$, $d\ge 0$, be the
blow@-up of $V_d$ along the proper transforms of the $x^d_r$'s and
$y^d_r$'s. Blowing up along the subspaces of dimension $d$ separates
the proper transforms of the subspaces of dimension $d+1$ containing
them, so at each stage the centers of the blow@-ups are necessarily
disjoint, and the blow@-ups can be performed in any order: in other
words, these varieties do not depend on the specific ordering given to
the $x$'s and $y$'s (in each dimension). Since $G$ is finite, this
construction stops at some stage, and we let $V_G$ be the resulting
variety. Of course $V_G$ depends on the dimension $n$ of the initial
projective space $\P^n$; however, in most of the paper this will not
play a r\^ole.

In $V_G$ we single out several natural divisor classes: the pull@-back
$H_0$ of the hyperplane class from $\P^n$; the pull@-backs $E^d_r$ of
the exceptional divisors arising by blowing up along $x^d_r$; the
pull@-backs $F^d_r$ of the exceptional divisors arising by blowing up
along $y^d_r$; and the classes $H^d_r$ of the proper transforms of
the general hyperplanes containing $x^d_r$. We define a divisor class
$S(t)$ as follows: let $R$ be the dimension of the subspace $x^R$
spanned by {\it all\/} the $x^0_r$ ($R+1$ equals then the number of
edges in a spanning forest of $G$; equivalently, the number of
vertices of the graph minus the number of its connected
components---cf.~\cite{W1}, 6.1.2); then set
$$S(t)=t^{R+1} H_0-\sum_{d,r} t^{R-d} E^d_r\quad.$$

\remark{Remark} Notice that the $F's$ are not used in this definition:
in fact, most computations in the following can be performed `modulo
$F$' (that is, modulo combinations of $F^d_r$'s). A construction could
be concocted without introducing the auxiliary subspaces $y^d_r$ and
the corresponding $F$'s, and still obtaining many of the results of
the paper. We have chosen this alternative path because the
construction as presented here is more natural in that it is
independent of the ordering of the subspaces, and moreover blowing-up
along the $y$'s makes the $H^d_r$'s generated by global sections (in
fact, this amounts to resolving at one time all maps defined in terms
of line bundles corresponding to nonnegative combinations of the
$H^d_r$'s). Is there an equally natural construction that does not
invoke the use of these `auxiliary' subspaces and
divisors?\endremark\vskip 12pt

The following is the prototype of the results in the paper. We defer
more general statements (and proofs) to later sections. Observe that
$$H_0 \text{, the } H^d_r \text {, and the } F^d_r$$
give a basis of the $\Pic$ of $V_G$. Now by Poincar\'e duality we can
find a class $\gamma\in A_1(V_G)$ dual to $H_0$ with respect to this
basis: that is, such that
$$H_0\cdot\gamma = 1, \quad H^d_r\cdot\gamma = 0, \quad
F^d_r\cdot\gamma = 0\quad\text{for all $d,r$.}$$
In other words, given a divisor $D$ in $V_G$, $D\cdot\gamma$ picks the
coefficient of $H_0$ in the (unique) expression of $D$ in terms of
$H_0$, $H^d_r$'s, and the $F^d_r$'s.

\proclaim{Theorem 1.1} Let $c$ be the number of connected components of
$G$. Then the number of ways in which $G$ can be colored properly with
$m$ colors (that is, so that no two adjacent vertices are given the
same color) is given by the intersection product
$$m^c\,S(m)\cdot\gamma\quad.$$
\endproclaim

\proclaim{Corollary 1.2} $G$ can be colored properly with $m$ colors if
and only if $S(m)\cdot\gamma \ne 0$.\endproclaim

\example{Examples} (1) If $G$ has at least 1 edge, then
$S(1)=H_0-\sum_{d,r} E^d_r$ is, modulo $F$, the class of the
proper transform of the hyperplane containing all the $x^d_r$'s; so
(by definition of $\gamma$) $S(1)\cdot \gamma = 0$. If $G$ has no
edges, then $V_G=\P^n$, $S(1)=H_0$, and thus $S(1)\cdot \gamma=1$. The
corresponding facts about proper colorings are of course trivial.

(2) Let $G$ be the complete graph on 4 vertices. The six $x^0_r$ are
placed at the points of intersection of four general lines of a plane;
on each of these four lines $x^1_1,\dots,x^1_4$ lie three of the
$x^0_k$. There are three pairs of $x^0_k$'s not lying on the same one
line in this configuration; these pairs determine three more lines
$x^1_5,x^1_6,x^1_7$. Finally, there is one plane $x^2$ containing the
whole configuration. By using the definition of $\gamma$, we find
$$\gather
E^0_k\cdot\gamma = 1, \quad k=1,\dots,6;\\
E^1_r\cdot\gamma = -2, \quad r=1,\dots,4, \quad\text{and }
E^1_r\cdot\gamma = -1, \quad r=5,6,7;\\
E^2\cdot\gamma = 6\quad,
\endgather$$
so
$$\align
m\, S(m)&=m\,(m^3-6\cdot 1\, m^2-(-2\cdot 4-1\cdot
3)\,m-1\cdot 6)\\
&=m^4-6 m^3+11 m^2- 6m=m\,(m-1)\,(m-2)\,(m-3)
\endalign$$
as it should be: each vertex must be assigned a different color from
the palette.
\endexample

We can prove a stronger statement than Theorem 1.1, which exploits one
of the basic features of the construction: $V_G$ encodes at once the
combinatorial information of $G$ and of all its {\it contractions\/.}
Each $x^d_r$ corresponds to a choice of edges of the original graph;
let $G^d_r$ be the graph obtained from $G$ by contracting each edge in
this collection, and removing parallel edges that might be created in
the process (note: no loops arise by this operation). Also, let
$\gamma^d_r$ be the dual of $H^d_r$ in the above basis. Up to a power
of $m$, then, {\it $S(m)\cdot\gamma^d_r$ counts the proper
$m$-colorings of the contraction $G^d_r$\/.} (This will follow from
Theorem 2.3 in the more general setting of \S2).

In other words, denote by $\overline S(m)$ the divisor equivalent to
$S(m)$ modulo $F$ and in the span of $H_0,H^d_r$: then the above says
that

\noindent{\it $G$ and all its contractions can be colored properly with
$m$ colors if and only if $\overline S(m)$ is in the interior of the cone
generated by $H_0, H^d_r$ in $\Pic V_G$\/.}

\noindent For example, the four@-color@-theorem (\cite{A-K}) says that
if $G$ is a planar graph, then $\overline S(4)$ is in the interior of
the cone generated by $H_0, H^d_r$ (since all contractions of a planar
graph are planar).

We end the section by remarking that in the case we have considered
here (that is, varieties arising from graphs), the $\overline S(m)$,
$m>0$, turn out to be all generated by global sections (see
Proposition 2.4): indeed, the $H^d_r$'s are, and, by the above
results, $\overline S(m)$ is a nonnegative combination of the $H_0$
and the $H^d_r$'s in the graph case. This does not seem obvious a
priori, for it is {\it not\/} true for the analogous construction for
matroids examined in the next section (we will find there a class of
matroids for which this holds, cf.~Proposition 2.5). In the graph
case, it follows that for positive $m$ there always is a hypersurface
in $\P^n$ generically smooth along the maximal $x^R$, with
multiplicity $m$ along the $x^{R-1}_r$'s, multiplicity $m^2$ along the
$x^{R-2}_r$'s, \dots, multiplicity $m^R$ at the $x^0_r$'s and degree
$m^{R+1}$: simply take general hyperplanes containing the $x^d_r$'s as
dictated by the expression of $\overline S(m)$ in terms of $H_0$ and
the $H^d_r$'s. The class $\overline S(m)$ is then the class of the
proper transform of such a hypersurface.

Conversely, we may view the above as a recipe to compute the chromatic
polynomial of a graph: given the collection of $x^d_r$'s obtained as
above, construct a hypersurface by taking enough general hyperplanes
containing each $x^d_r$ to satisfy the above multiplicity prescription
(multiplicity 1 along the maximal subspace $x^R$, $t$ along
codimension 1 subspaces, $t^2$ along codimension 2, etc.). By the
above, this will always be possible: the number needed at $x^d_r$ is
$S(t)\cdot \gamma^d_r\ge 0$; and the number of hyperplanes not
containing any of the $x^d_r$'s, needed to get a hypersurface of
degree $t^{R+1}$, multiplied by $t$ to a power equal to the number of
connected components of $G$, will give the value at $t$ of the
chromatic polynomial of $G$ (this is of course nothing but ``M\"obius
inversion'' at work).\vskip 12pt


\heading \S 2. Matroid varieties\endheading

In section 1 we gave the standard embedding in a projective space of
the `cycle matroid' associated with the graph $G$, and constructed a
variety $V_G$ from this data. The construction can be performed for
the lattice $\Cal L=\Cal L(\Cal C)$ of subspaces spanned by any finite
collection $\Cal C$ of points in $\P^n$.  $\Cal L$ is (partially)
ordered by inclusion; 0~will be the empty set (the minimum of the
lattice), 1 the maximal subspace, spanned by all points; we require
this to have codimension at least 2 in $\P^n$. We denote elements of
$\Cal L$ by letters $x,y,z,\dots$, by $\le$ the ordering in $\Cal L$,
and by $\vee,\wedge$ resp.~the join and meet in the lattice. The
`rank' $r(x)$ of $x\in\Cal L, x\ne 0,$ is one plus its dimension as a
subspace of $\P^n$: so the points of $\Cal C$ are the rank@-1 elements
of $\Cal L$. The rank of $0=\emptyset$ is 0; the `rank of $\Cal L$' is
$r(\Cal L)=r(1)$.

Now $V_{\Cal L}$ is constructed as in section 1. First we close the
family $\Cal L$ of subspaces of $\P^n$ with respect to intersection:
let $\Cal M$ be the family of subspaces $\notin \Cal L$ obtained by
intersecting collections of elements of $\Cal L$; we extend rank and
ordering to elements of $\Cal M$. Next, $V_{\Cal L}$ is obtained from
$\P^n$ by blowing up the (proper transforms of the) $x\ne 0$ in $\Cal
L$ and $\Cal M$ in order of increasing dimension; again we observe
that since $\Cal L\cup\Cal M$ is closed under intersections,
blowing@-up all $x$ of rank $r$ separates the proper transforms of the
subspaces of rank $r+1$, hence the construction is independent of the
specific order in which the blow@-ups are executed (within each rank).

We note that $V_G=V_{\Cal L}$ if $\Cal L$ corresponds to $G$ as in
section 1. Keeping the same style of notations as in \S 1, we let
$H_x$ be the class of the proper transform of the general hyperplane
containing $x$ (so the pull@-back of the hyperplane class is $H_0$),
we let $E_x$ be the pull@-back of the exceptional divisor over $x\in
\Cal L, x\ne 0$, and $F_x$ be the pull@-back of the exceptional
divisor over $x\in \Cal M$. For $x\in\Cal L$, $\gamma_x$ is a 1@-class
such that $\gamma_x\cdot H_x=1,\gamma_x\cdot H_y=0$ for all $y\in\Cal
L, y\ne x$, and $\gamma_x\cdot F_z=0$ for all $z\in\Cal M$. $S(t)$ is
the class
$$S(t)= t^{r(1)} H_0-\sum_{x\in\Cal L,x\ne 0} t^{r(1)-r(x)} E_x$$

\noindent(as in \S 1, we will soon introduce a class $\overline S(t)$
equivalent to $S(t)$ `modulo $F$' but somewhat better behaved.)

\subheading{\S 2.1. Compatibilities with contractions, deletions,
etc\/} We will now show how the construction behaves with respect to
three basic matroid operations. All the results in \S3 will
essentially follow from a closer look to the compatibilities sketched
below; a detailed analysis of the functorial properties of the
construction is well beyond the scope of this note. For the hasty
reader: only contractions will be used in the rest of this
section.\vskip 12pt

{\it Contractions\/.} The variety $V_{\Cal L}$ contains a `compatible'
copy of $V_{\Cal L/x}=V_{[x,1]}$ for each $x\in L$. More precisely:
the fiber of the exceptional divisor obtained when blowing@-up along
$x\in\Cal L$ is a projective space $\P^{n-r(x)}$, met by all and only
the $z\ge x$ in $\Cal L$. The lattice of subspaces these form in this
projective space is the interval $[x,1]$, isomorphic to the `geometric
contraction' $\Cal L/x$ of $\Cal L$ by $x$ (\cite{W1}, p.~141). In
terms of graphs, this is the contraction determined by a choice of a
collection of edges, as described in \S 1. Now the blow@-up process is
compatible with restriction to this $\P^{n-r(x)}$: the general fiber
of $E_x$ (that is, the proper transform of $\P^{n-r(x)}$ in $V_{\Cal
L}$) is the blow@-up of $\P^{n-r(x)}$ along its intersection with the
$z\in \Cal L\cup\Cal M$, $z\ge x$, that is nothing but a copy of
$V_{\Cal L/x}$. Further, all expected compatibilities among the
definitions of the relevant classes hold; for example, the class
$\gamma_x$ in $V_{\Cal L}$ is the push@-forward of the class
$\gamma_0$ in $V_{\Cal L/x}$, etc. Typically, anything proved about
$\Cal L$ by means of $V_{\Cal L}$ will automatically restrict to a
statement about all its contractions.\vskip 6pt

{\it Modular elements\/.} At the same time, $V_{\Cal L}$ also contains
a copy of $V_{[0,x]}$ (where $[0,x]$ denotes the lattice of elements
$z\in\Cal L$ such that $0\le z\le x$), provided that $x$ be {\it
modular\/.} An element $x\in\Cal L$ is `modular' if $x\wedge z=x\cap
z$ for all $z$ in $\Cal L$ (where $\wedge$ denotes the meet in the
lattice, while $\cap$ denotes intersection in $\P^n$); for example,
all rank@-1 elements of $\Cal L$ are modular. Now consider any
subspace $\P_x$ of $\P^n$, of dimension $>r(x)$ and intersecting
$1\in\Cal L$ precisely along $x$; then
\proclaim{Claim 2.1} If $x$ is modular, then the
proper transform of $\P_x$ in $V_{\Cal L}$ is isomorphic to a variety
$V_{[0,x]}$.\endproclaim
\demo{Proof} $\P_x$ contains a copy of $[0,x]$. Let $\Cal M_x$ denote
for a moment the set of subspaces defined when constructing
$V_{[0,x]}$ (that is, all $y\cap z\notin [0,x]$, where $y,z\in
[0,x]$). Then it is easily checked that modularity implies
$[0,x]=\{z\cap \P_x, z\in\Cal L\}$ and $\Cal M_x=\{z\cap \P_x,
z\in\Cal M\}$. Taking the proper transform of $\P_x$ amounts then to
performing precisely the same sequence of blow-ups producing
$V_{[0,x]}$ as dictated by the construction.\qed\enddemo

{\it Deletions.\/} The construction is also compatible with
substructures. Let $\Cal C'$ be a subset of the set of rank-1 elements
of $\Cal L$ (that is, of the original set $\Cal C$ of points in $\P^n$
generating~$\Cal L$); these generate a sublattice $\Cal L(\Cal C')$ of
$\Cal L$, a `deletion' of $\Cal L$. Then there is a map $V_{\Cal L}
@>>> V_{\Cal L(\Cal C')}$: this follows from the universal property of
blow@-ups, once we observe that the inverse image of all subspaces
generated by elements of $\Cal C'$ (and all their intersections) are
Cartier divisors in $V_{\Cal L}$. For example, for $\Cal
C'=\emptyset$, the resulting map $V_{\Cal L}@>>> V_{\Cal
L(\emptyset)}=\P^n$ is simply the sequence of blow@-ups
defining~$V_{\Cal L}$.\vskip 6pt

{\it Nesting\/.} Finally, we observe that we get a variety
$V^n=V_{\Cal L}$ by blowing up $\P^n$ as above, for each $n>r(1)$;
most results of the paper do not depend on the specific choice of $n$.
These different varieties are nested into each others like Russian
dolls: for all $n>r(1)$, $V^n$ can be embedded as a divisor of class
$H_1$ in $V^{n+1}$. Indeed, the proper transform of any $\P^n$
containing $1\in \Cal L$ in $\P^{n+1}$ is a copy of $V^n$: this is
Claim 2.1 for $x=1$ (1~is always modular!).

\subheading{\S 2.2. The characteristic polynomial} Now for a bit of
well known and beautiful combinatorics, and its translation into the
intersection ring of $V_{\Cal L}$. Recall (\cite{W2}, Chapter 7) that
the `M\"obius function' of a lattice $\Cal L$ is the function
$\mu_{\Cal L}: \Cal L\times \Cal L @>>> \Z$ satisfying
$$\sum_{x\le y\le z}\mu_{\Cal L}(x,y)=\left\{\gathered\text{0 if $x\ne
z$} \\\text{1 if $x=z$}\endgathered\right.\quad\text{if $x\le z$},\quad
\text{$\mu_{\Cal L}(x,z)=0$ if $x\not\le z$}$$
We will write $\mu$ for $\mu_{\Cal L}$ if no ambiguity is feared. The
`characteristic polynomial' of $\Cal L$ is the polynomial
$$p(\Cal L,t)=\sum_{x\in \Cal L}\mu(0, x)\,t^{r(1)-r(x)}$$
Now the key observation is the following (see for example \cite{W2},
\S 7.5): the number of proper colorings of a graph $G$ with $t$ colors
(that is, the `chromatic polynomial' of $G$) is given by
$$t^c\, p(\Cal L,t)\quad,$$
where $c$ is the number of connected components of $G$ and $\Cal L$ is
the lattice determined by $G$. So Theorem 1.1 will be proved once we
show that for any matroid in $\P^n$ as above:
\proclaim{Theorem} $p(\Cal L,t)=S(t)\cdot\gamma_0\quad.$\endproclaim

In turn, given the definition of $S(t)$, this is proved once we
observe that $\mu(0,0)=1=H_0\cdot\gamma_0$, and show that
$\mu(0,z)=-E_z\cdot \gamma_0$ for $z\in\Cal L, z\ne 0$. In fact:
\proclaim{Lemma 2.2} $E_z \cdot\gamma_x=-\mu(x,z)$ for all $z\in\Cal
L, z\ne 0$.\endproclaim
\demo{Proof} First we observe that by restricting to the general fiber
of $E_x$ we may assume $x=0$ (by compatibility with contractions).
So we just have to show $\mu(0,z)=-E_z\cdot \gamma_0$ for
$z\ne 0$. By definition of the M\"obius function, this amounts to
showing
$$\mu(0,0)+\sum_{0< y\le z}(-E_y\cdot\gamma_0)=0$$
whenever $z\ne 0$. But observe that the construction gives
$$H_z = H_0 - \sum \Sb y\in\Cal L \\ 0< y\le z \endSb E_y - \sum \Sb
x\in\Cal M \\ x< z \endSb F_x \quad,$$
so that
$$\align
\mu(0,0)+\sum_{0< y\le z}(-E_y\cdot\gamma_0)&= 1+\sum_{0< y\le
z}(-E_y\cdot\gamma_0)\\
&=(H_0-\sum_{0< y\le z} E_y -\sum_{x<z} F_x)\cdot \gamma_0\\
&=H_z\cdot\gamma_0 = 0
\endalign$$
by definition of $\gamma_0$.\qed\enddemo

As pointed out, this lemma implies the theorem above, and this in turn
implies Theorem 1.1. There is a substantial advantage, however, in
giving a more comprehensive statement dealing with all contractions of
$\Cal L$ at once. For this, let $\overline S(t)$ denote the divisor
equivalent to $S(t)$ modulo $F$ and in the span of the $H_x$'s. Note
that $S(t)$ and $\overline S(t)$ have the same intersection numbers
against any combination of the $\gamma_x$, $x\in\Cal L$.

\proclaim{Theorem 2.3} Denoting by $\Cal L/x\cong [x,1]$ the
sublattice of $\Cal L$ consisting of all $z\in\Cal L$ such that $x\le
z$:
$$S(t)\cdot\gamma_x=p(\Cal L/x,t)$$
for all $x\in\Cal L$. In other words,
$$\overline S(t)=\sum_{x\in\Cal L} p(\Cal L/x,t)\, H_x\quad.$$
\endproclaim
\demo{Proof}
$$\align
S(t)\cdot\gamma_x&=t^{r(1)} H_0\cdot\gamma_x-\sum_{y\in\Cal
L,x\ne 0} t^{r(1)-r(y)} E_y\cdot\gamma_x\\
&=\sum_{y\in\Cal L} t^{r(1)-r(y)}\mu(x,y)\qquad\text{by Lemma 2.2}\\
&=\sum_{y\ge x} t^{r(1)-r(y)}\mu(x,y)\\
&=p([x,1],t)=p(\Cal L/x,t)\quad.\qed
\endalign$$
\enddemo

Theorem 2.3 implies the extension of Theorem 1.1 discussed in section 1.

$S(t)$ is easier to define, while $\overline S(t)$ is better behaved
in some respects. For example, $\overline S(m)$ is automatically
globally generated for $m>0$ in the graph case (as mentioned in \S1),
because:

\proclaim{Proposition 2.4} Non@-negative linear combinations of the
$H_x$'s are generated by global sections.\endproclaim

\demo{Proof} We only need to show that each $H_x$ is generated by
global sections. Now $H_0$ clearly is, since it `already' is in
$\P^n$; for $x\ne 0$, observe that any $x\in\Cal L$ is the
intersection of $n+1-r(x)$ general hyperplanes containing it. In the
construction, every center of blow@-up is either included in the
proper transform of $x$, or it is disjoint from it (note: this would
not necessarily be the case if we didn't blow@-up along the elements
of $\Cal M$ as well!). It follows that the proper transforms of the
hyperplanes still intersect exactly along the proper transform of $x$
after each blow@-up, and get separated when $x$ itself is blown up.
They give then $n+1-r(x)$ sections of $H_x$ generating it
globally.\qed\enddemo

\remark{Remark} What was shown in this proof was in fact that
$n+1-r(x)$ general representatives of $H_x$ have empty intersection in
$V_{\Cal L}$.\endremark
In the graph case, the coefficients of $H_x$ in $\overline S(t)$ are
(up to powers of $t$) chromatic polynomials, thus nonnegative at
positive integers: so $\overline S(m)$ is in the cone generated by the
$H_*$ in $\Pic V_G$ for all positive $m$, and is globally generated.

This does not seem at all obvious a priori, say from the definition of
$S(t)$; in fact, it is {\it not\/} true for arbitrary matroids! For
example, consider the matroid $L_4$ generated by four collinear
points: if $\overline S(2)$ were generated by global sections, then
(at least in char.~0) by Bertini there would be a nonsingular
irreducible hypersurface of class $\overline S(2)$ in $V_{L_4}$; this
would map down to $\P^n$ to a hypersurface of degree 4, generically
smooth along a line, and having multiplicity 2 at (at least) 4 points
on this line.  This cannot be: the general plane section of this
hypersurface would be a plane quartic curve containing a line, whose
residual cubic meets the line at four distinct isolated points. Thus
$\overline S(2)$ is not generated by global sections in general.

It would be interesting to find a characterization of planar graphical
matroids in terms of properties of the divisors $\overline S(m)$. A more
ambitious goal would be to find for each given matroid $M$ an
algebro@-geometric property of $V_{\Cal L}$ that can signal whether $\Cal
L$ is the lattice of a matroid {\it none of whose minors is isomorphic
to $M$.\/} Such a tool would allow us to mirror the characterization
of classes of matroids in terms of `excluded minors' (see pp.~146--7
in \cite{W1}); in particular a characterization of varieties arising
from planar graphical matroids would follow.

The only result of this sort that we know is the following. Following
the common terminology, we denote by $L_4$ the `four point line' of
the above example, and by $F_7$ the `seven point plane' (that is, the
matroid defined by the projective plane over the 2@-element field).

\proclaim{Proposition 2.5} Let $\Cal L$ be the lattice corresponding
to a given matroid $M$, and $\overline S(m)$ the divisor on $V_{\Cal
L}$ defined as above. Then the following are equivalent:\roster
\item $\overline S(2)$ and $\overline S(3)$ are in the cone generated
by the $H_x$, $x\in\Cal L$;
\item $M$ has no minor isomorphic to $L_4$ or $F_7$;
\item All $\overline S(m)$, $m>0$, are in the cone generated by the
$H_x$, $x\in\Cal L$.
\endroster\endproclaim

\remark{Remark} This amounts to saying that the class defined in (2)
is precisely the class of matroids whose contractions all have
characteristic polynomials which are non@-negative at each positive
integer. This must be a well@-known characterization in combinatorics,
but we could not trace it in the literature; we apologize for the
missing reference and provide the following straightforward (and
hopefully correct) argument.\endremark

\demo{Proof} (3) $\implies$ (1) is trivial.

(1) $\implies$ (2): if $M$ has a minor isomorphic to $L_4$, then by
the `scum theorem' (Prop.~7.4.11 in \cite{W1}) $L_4$ is obtained from
$M$ by a contraction $M/I$ followed by a sequence of deletions:
$L_4=M/I-e_1-\dots-e_r$. Now $p(L_4,m)=m^2-4 m+3$, so $p(\Cal
M/I-e_1-\dots-e_r,2)=p(L_4,2)=-1$; we claim that this implies some
contraction of $M$ has negative characteristic polynomial at 2.
Indeed, by \cite{W2}, Theorem 7.2.4,
$$\multline
p(M/I-e_1-\dots-e_r,2)=p(M/I-e_1-\dots-e_{r-1},2)\\
+ p(M/(I\vee e_r)-e_1-\dots-e_{r-1},2)
\endmultline$$
if $e_r$ is not an isthmus in $M/I-e_1-\dots-e_{r-1}$, and
$$p(M/I-e_1-\dots-e_r,2)=p(M/I-e_1-\dots-e_{r-1},2)$$
if $e_r$ is an isthmus in  $M/I-e_1-\dots-e_{r-1}$. In either case,
the polynomial is necessarily negative at 2 for a contraction of $M$
followed by fewer deletions: the claim follows. Finally, the
coefficients in the expression of $\overline S(2)$ in terms of the
$H_*$ are precisely the values of the characteristic polynomials of
the (geometric) contractions of $M$ (by Theorem 2.3), so we can
conclude that $\overline S(2)$ is not in the cone generated by the
$H_*$. The argument for $F_7$ is entirely similar, given that
$p(F_7,m)=m^3-7 m^2+14 m-8$ is negative for $m=3$.

(2) $\implies$ (3): the class defined in (2) is closed under
contractions, so we just need to show that the characteristic
polynomial of any matroid in it is nonnegative at positive integers.
By a result of Seymour (cf.~\cite{W1}, p.~147), the class is in fact
the class of `direct sums and 2@-sums of regular matroids and copies
of $F^*_7$'. Now observe that $p(F^*_7,m)=m^4-7m^3+21 m^2-28 m+13$ is
$\ge 0$ for all integer $m>0$; also, regular matroids have nonnegative
characteristic polynomial because of a result of Crapo (Theorem III in
\cite{C}: the value of the polynomial at $m$ is the number of
`$H$@-coboundaries with kernel 0', for $H$ a group of order $m$).
Next, nonnegativity is preserved by direct sums by Theorem 7.2.4 (ii)
in \cite{W2}; so we just have to show it is preserved under 2@-sums.
Now the 2@-sum of two matroids $M_1,M_2$ is obtained from their
parallel connection by deletion of the base point: in the notation of
\cite{W1}, p.~180
$$S_2(M_1,M_2)=P(M_1,M_2)-p\quad,$$
where $p$ is not an isthmus of either $M_1$ or $M_2$. It follows that
$p$ is not an isthmus of $P(M_1,M_2)$, so applying 7.2.4 (i) from
\cite{W2}, $7.6.7_P$ from \cite{W1}, 7.2.9 and 7.2.4 (ii) from
\cite{W2} we get
$$\align p(S_2(M_1,M_2),m)&=p(P(M_1,M_2),m)+p(P(M_1,M_2)/p,m)\\
&=p(P(M_1,M_2),m)+p(M_1/p\oplus M_2/p,m)\\
&=\frac{p(M_1,m)\,p(M_2,m)}{m-1}+p(M_1/p,m)\,p(M_2/p,m)\quad:
\endalign$$
each summand on the right is non@-negative, so we are done.\qed\enddemo

All matroids representable over any field, and in particular all
graphical matroids, belong to the class defined in this proposition;
however, for such matroids one can prove (3) more directly, cf.~the
discussion following Proposition 2.4. For all matroids satisfying (3),
the line bundles corresponding to $\overline S(m)$ are globally
generated, so they define maps from the variety of the matroid to a
projective space. We feel that studying these maps would be quite
fruitful; we will obtain a simple result about the degree of the image
of such maps in \S 4.

Of course a characterization of planar graphs in a fashion similar to
Proposition~2.5 would be desirable.\vskip 12pt


\heading \S 3. Characteristic polynomial basics, Crapo's invariant: a
geometric viewpoint\endheading

In this section we run through basic material concerning
characteristic polynomials, illustrating it in the context of the
construction introduced in \S2. The reader is encouraged to compare
the `geometric' proofs given here with more standard combinatorial
arguments, as presented for example in Chapter 7 of \cite{W2}.

The general strategy is the following: in a given situation, write the
most fundamental relation suggested by the geometry; then applying the
results in \S 2 will yield an equally fundamental combinatorial
statement. As an appetizer, the following is the simplest possible
example of such an argument:

\proclaim{Proposition 3.1} With notations as in \S2\/, $\sum_{x\in\Cal
L}\gamma_x$ equals the class of the pull-back $\ell$ of a line from
$\P^n$.\endproclaim
\demo{Proof} Dot both classes against all divisors.\qed\enddemo

And here is the translation into combinatorics:

\proclaim{Corollary 3.2} $\sum_{x\in\Cal L} p(\Cal L/x,t) = t^{r(1)}$
\endproclaim
\demo{Proof} By Theorem 2.3 and Proposition 3.1, the left@-hand@-side
is $S(t)\cdot\sum_{x\in\Cal L}\gamma_x=S(t)\cdot\ell$. But the
pull@-back of a line vanishes against all exceptional divisors, so
$S(t)\cdot\ell = t^{r(1)} H_0\cdot\ell=t^{r(1)}$.\qed\enddemo

The other examples in this section are a little more complex, but
motivated by the same simple geometric intuition.

\subheading{\S 3.1. Deletion-contraction rule} Let $e\in\Cal C$ be a
rank@-1 element in $\Cal L$---that is, one of the points in the set
used to generate the subspaces in $\Cal L$. Denote by $\Cal L-e$ the
lattice of subspaces spanned by the {\it other\/} points ($\Cal L-e$
is a `deletion' of $\Cal L$). We observed in \S2.1 that the universal
property of blow@-ups gives then a map
$$\alpha: V_{\Cal L} @>>> V_{\Cal L-e}$$
compatible with the blow@-up maps from the matroid varieties to
$\P^n$. In particular, this map is proper, birational and onto. We use
notations as in \S2, and append a $'$ to denote objects in $V_{\Cal
L-e}$: so e.g., $H_0'$ is the pull-back of the hyperplane class to
$V_{\Cal L-e}$ (and it follows $\alpha^*(H_0')=H_0$), etc.

\proclaim{Proposition 3.3} $\alpha^*(\gamma_0')=\gamma_0+\gamma_e$
\endproclaim
\demo{Proof} It is clear that the class vanishes against `$F$
divisors'; we have to show $H_x\cdot \alpha^*(\gamma_0)=0$ if $x\ne
0,e$, and $=1$ otherwise. Now any $x\in\Cal L$, $x\ne 0, e$, contains
a maximal $x'\in\Cal L-e$, $x'\ne 0$; the reader will then check that
$\alpha_*(H_x)=H'_{x'}$. Since $\alpha$ is birational, and using the
projection formula, $H_x\cdot\alpha^*(\gamma_0')=H'_{x'}\cdot
\gamma_0'=0$ since $x'\ne 0$. By the same token,
$\alpha_*(H_e)=\alpha_*(H_0)=H_0'$, from which
$H_e\cdot\alpha^*(\gamma_0')=H_0\cdot\alpha^*(\gamma_0')=1$.\qed\enddemo

Proposition 3.3 `stands behind' the deletion@-contraction rule for the
characteristic polynomial (Theorem 7.2.4(i) in \cite{W2}), curiously
regardless of $e$ being or not an {\it isthmus\/} of $\Cal L$ (a
rank@-1 element $e$ of $\Cal L$ is an `isthmus' if the rank of $\Cal
L$ is strictly larger than the rank of $\Cal L-e$). More precisely:

\proclaim{Corollary 3.4(a)} If $e$ is not an isthmus, then $p(\Cal
L,t)=p(\Cal L-e,t)-p(\Cal L/e,t)$.\endproclaim
\demo{Proof} If $e$ is not an isthmus, then $r(\Cal L)=r(\Cal L-e)$,
and it follows that $\alpha_*(S(t))=S(t)'$ by definition. By Theorem
2.3 and using the projection formula:
$$\align
p(\Cal L-e,t)&=S(t)'\cdot\gamma_0'=\alpha_*(S(t))\cdot\gamma_0'
=S(t)\cdot\alpha^*(\gamma_0')\\
&=S(t)\cdot(\gamma_0+\gamma_e)\qquad\text{by the proposition}\\
&=p(\Cal L,t)+p(\Cal L/e,t)
\endalign$$
again by Theorem 2.3.\qed\enddemo

\proclaim{Corollary 3.4(b)} If $e$ is an isthmus, then $p(\Cal
L,t)=(t-1)\, p(\Cal L-e,t)$.\endproclaim
\demo{Proof} If $e$ is an isthmus, then $r(\Cal L)=r(\Cal L-e)+1$.
 From this it follows that $\alpha_*(S(t))=t \,S(t)'$, so
$$\align
t\, p(\Cal L-e,t)&=t \, S(t)'\cdot\gamma_0'=\alpha_*(S(t))
\cdot\gamma_0'\\
&=S(t)\cdot(\gamma_0+\gamma_e)\qquad\text{arguing as above}\\
&=p(\Cal L,t)+p(\Cal L/e,t)\\
&=p(\Cal L,t)+p(\Cal L-e,t)
\endalign$$
since $\Cal L/e=\Cal L-e$ if $e$ is an isthmus. The statement
follows.\qed\enddemo

\subheading{\S 3.2. Stanley's modular factorization theorem} If $\Cal
L$ is the product $\Cal L_1\times\Cal L_2$ of two lattices, we
could argue as above and prove the multiplicativity of the
characteristic polynomial under direct sums, by studying the map
$V_{\Cal L} @>>> V_{\Cal L_1}$. However, as pointed out in \cite{W2},
p.~122, this is a particular case of a more general factorization
result (\cite{S}, Theorem 2); so we present the latter.

Recall from \S 2 that we have an injection $i: V_{[0,x]}
\hookrightarrow V_{\Cal L}$ whenever $x$ is a modular element of $\Cal
L$. Again we use notations as in \S 2, appending a $''$ to denote
objects of $V_{[0,x]}$.
\proclaim{Proposition 3.5}If $x$ is modular, and with notations as
above:\roster
\item $i^*(H_z)=H_{x\wedge z}''$;
\item $i^*(E_z)=E_z''$ if $z\le x$, 0 otherwise;
\item $i^*(F_z)=F_z''$ if $z\le x$, 0 otherwise;
\item $i^*(S(t))=t^{r(1)-r(x)}\,S(t)''\quad\text{
and }\quad i^*(\overline S(t))=t^{r(1)-r(x)}\,\overline S(t)''$
\endroster\endproclaim
\demo{Proof} (2) and (3) follow from a chase of the diagram of
blow@-ups producing the two varieties. For example, if $z\nleq x$ then
$x$ and $z$ are separated when blowing@-up along $x\cap z$; the
proper@-transform of $\P_x$ (that is, $V_{[0,x]}$ by Claim 2.1) is
then disjoint from the exceptional divisor above $z$, and the
corresponding pull@-back must vanish.

(1) follows from (2) and (3). The first part of (4) follows from the
definitions of $S(t)$, $S(t)''$ and from (1) and (2). The second part
of (4) follows from the first, by killing $F$ terms on both
sides.\qed\enddemo

\proclaim{Corollary 3.6} For all modular $x\in\Cal L$ and all
$y\in[0,x]$
$$\sum_{z\in\Cal L, z\wedge x=y} p(\Cal L/z, t)=t^{r(1)-r(x)}\,
p([y,x],t)$$
\endproclaim
\demo{Proof} Using Theorem 2.3 to write out the second part of (4)
from the proposition:
$$\align
t^{r(1)-r(x)}\sum_{0\le y\le x} p([0,x]/y,t)\,H_y''&=i^*(
\sum_{z\in\Cal L}p(\Cal L/z,t) \,H_z)\\
&=\sum_{z\in\Cal L}p(\Cal L/z,t) \,H_{z\wedge x}''\qquad\text{by (1)
above}\\
&=\sum_{0\le y\le x}\left(\shave{\sum_{z\in\Cal L, z\wedge x = y}}
p(\Cal L/z,t)\right) H_y''\quad.
\endalign$$
The statement follows by dotting with $\gamma_y''$ and observing
$[0,x]/y=[y,x]$.\qed\enddemo
Setting $y=0$ in the statement and isolating $p(\Cal L/0,t)=p(\Cal
L,t)$ gives
$$p(\Cal L,t)=t^{r(1)-r(x)}\,p([0,x],t)-\sum \Sb z\in\Cal L,z\ne 0\\
z\wedge x=0\endSb p(\Cal L/z,t)$$

\proclaim{Corollary 3.7} (Modular factorization theorem) If $x$ is a
modular element of $\Cal L$, then
$$p(\Cal L,t)=p([0,x],t)\sum_{y\in\Cal L, y\wedge x=0}\mu(0,y)\,
t^{r(1)-r(x)-r(y)}$$
\endproclaim
\demo{Proof} By induction on the rank of $\Cal L$. The statement is
clear if the rank of $\Cal L$ equals $r(x)$ (because this forces
$x=1$). If $x$ is modular in $\Cal L$ and $z\wedge x=0$, then $z\vee
x$ is modular in $[z,1]=\Cal L/z$; and $r(\Cal L/z)<r(\Cal L)$ if
$z\ne 0$, so we may assume the statement for $\Cal L/z$ in this case.
Doing so in the formula preceding the statement of this corollary
gives the induction step.\qed\enddemo

\subheading{\S 3.3. Crapo's beta invariant} Writing down an expression
for the canonical divisor $\omega_{\Cal L}$ of a matroid variety
$V_{\Cal L}$ is an elementary exercise:
$$\omega_{\Cal L}=-(n+1)\,H_0+ \sum_{x\in\Cal L,x\ne
0}(n-r(x))\,E_x+\sum_{y\in \Cal M}(n-r(y))\,F_y\quad,$$
where $n$ denotes as usual the dimension of $V_{\Cal L}$. On the other
hand, an important invariant of a matroid is its {\it beta invariant}
$$\beta(\Cal L)=(-1)^{r(1)-1}\frac d{dt}\,p(\Cal L,1)$$
(our source is \S 7.3 in \cite{W2}). The beta invariant contains a
surprising amount of information: for example, it vanishes precisely
if the matroid is a direct sum (or it is trivial). Now it turns out
that the beta invariant of $\Cal L$ is intimately related to the
canonical divisor of $V_{\Cal L}$---thus its relevance is clear from
an algebro@-geometric perspective.
\proclaim{Proposition 3.8} Assume $\Cal L\ne 0$. Then
$$\beta(\Cal L)=(-1)^{r(1)}(1+\omega_{\Cal L}\cdot\gamma_0)$$
\endproclaim
We will see in a moment (Proposition 3.10) that knowing the
exceptional divisor of $V_{\Cal L}$ is in fact equivalent (modulo $F$)
to knowing the beta  invariant of all contractions of $\Cal L$.

\demo{Proof} $p(\Cal L,t)=S(t)\cdot\gamma_0$ (Theorem 2.3), so
$$\align
&(1+\omega_{\Cal L}\cdot\gamma_0)-(-1)^{r(1)}\beta(\Cal L)= 1+\left(
\omega_{\Cal L}+\frac{dS}{dt}(1)\right)\cdot\gamma_0\\
&=1+\left(\shave{-(n+1)\,H_0+\sum_{x\ne 0}(n-r(x))E_x+r(1)\,H_0-
\sum_{x\ne 0}(r(1)-r(x))\,E_x}\right)\cdot\gamma_0\\
&=1-\left(\shave{H_0+(n-r(1))\,\left(\shave{H_0-\sum_{x\ne 0}E_x}
\right)}\right)\cdot\gamma_0\\
&=1-(H_0+(n-r(1))\,H_1)\cdot\gamma_0\qquad\qquad\text{(modulo $F$)}\\
&=0
\endalign$$
as needed.\qed\enddemo

Notice that the canonical divisor depends on the dimension $n$ of
$V_{\Cal L}$; as the proposition shows, its intersection with
$\gamma_0$ does not (if $\Cal L\ne 0$). The reason is that each
variety is embedded in the next as a divisor of class $H_1$: so their
canonical divisors differ by multiples of $H_1$ by adjunction, and
their difference is not detected by $\gamma_0$ by definition of the
latter.

The excluded case ($\Cal L\ne 0$) and the shape of the formula in
Proposition 3.8 reflect a little white noise in the definitions. We
can improve the situation by modifying the definition slightly in
order to make it independent of $n$ and fully compatible with
contractions. To this effect, define the `beta divisor' of $V_{\Cal
L}$ to be
$$\Til\omega_{\Cal L}=H_0+(\dim V_{\Cal L}-r(\Cal L))\, H_1+
\omega_{\Cal L}$$
The beta divisor is more natural with respect to the construction, in
the sense that it is compatible with the operations we have
encountered so far. More precisely, let $\alpha: V_{\Cal L}@>>>
V_{\Cal L-e}$ and $i: V_{[0,x]} \hookrightarrow V_{\Cal L}$ be as in
\S\S 3.1,2 (with the same $'$, $''$ notations), and view $V_{\Cal
L/x}$ as a subvariety of $V_{\Cal L}$ as usual (\S 2.1); then
\proclaim{Proposition 3.9} For $e\in\Cal L,r(e)=1$ and $x\in\Cal
L$:\roster
\item $\Til\omega_{\Cal L}|_{V_{\Cal L/x}}=\Til\omega_{\Cal L/x}$
\item $\alpha_*(\Til\omega_{\Cal L})=\Til\omega_{\Cal L-e}-(r(\Cal
L)-r(\Cal L-e))H_1'$
\item if $x$ is modular, $i^*(\Til\omega_{\Cal L})=
\Til\omega_{[0,x]}-(r(1)-r(x))H_1''$
\endroster\endproclaim
\demo{Proof} These are all immediate from the definition and the
adjunction formula. For example, let's check (3): $V_{[0,x]}$ is
embedded in $V_{\Cal L}$ as the proper transform of a space
intersecting $1\in\Cal L$ precisely along $x$; it follows that
$V_{[0,x]}$ is cut out by $\dim V_{\Cal L}-\dim V_{[0,x]}$
representatives of $H_x$. Its normal bundle has then first Chern class
$=(\dim V_{\Cal L}-\dim V_{[0,x]})H_x$, so the adjunction formula and
Proposition 3.5(1) give
$$i^*(\omega_{\Cal L})=\omega_{[0,x]}-i^*((\dim V_{\Cal L}-\dim
V_{[0,x]})H_x)=\omega_{[0,x]}-(\dim V_{\Cal L}-\dim V_{[0,x]})H_x''$$

Plugging this into the definition of the beta divisor gives (3).
\qed\enddemo

These compatibility properties of the beta divisor are in our view the
motor behind the basic properties of the beta invariant (e.g., 7.3.1,
7.3.2 in \cite{W2}). To support this viewpoint, we derive a few of
these in the remaining of this section. For a start, let's observe
explicitly that knowing the beta divisor (modulo $F$) is equivalent to
knowing the beta invariant of $\Cal L$ and of all its contractions.
Indeed:

\proclaim{Proposition 3.10} For all $x\in\Cal L$:
$$\Til \omega_{\Cal L}\cdot\gamma_x= (-1)^{r(1)-r(x)}\beta(\Cal L/x)$$
\endproclaim
\demo{Proof} For $x=0$ this follows at once from Proposition 3.8, or
by explicit computation if $\Cal L=0$; the general case reduces to
$x=0$ by compatibility with contractions (Proposition 3.9(1)).
\qed\enddemo

$(-1)^{r(1)-r(x)}\beta(\Cal L/x)$ is called the `signed beta
function', $B(x)$, in \cite{W2}, \S 7.3. Proposition 3.10 simply says
$$\Til\omega_{\Cal L}=\sum_{x\in\Cal L}B(x)\, H_x\qquad\text{ modulo
$F$.}$$

\proclaim{Corollary 3.11} $\dsize\sum_{x\in\Cal L, x\ge
y}B(x)=r(y)-r(1)$\endproclaim
\demo{Proof} Reduce to $y=0$ by replacing $\Cal L$ by $\Cal L/y$. Then
by the last proposition
$$\align
\sum_{x\in\Cal L}B(x)&=\Til\omega_{\Cal L}\cdot\sum_{x\in\Cal L}
\gamma_x\\
&=\Til\omega_{\Cal L}\cdot\ell\qquad\text{by Proposition 3.1}\\
&=1+(n-r(1))-(n+1)=-r(1)
\endalign$$
by the definition of $\Til\omega_{\Cal L}$.\qed\enddemo

The above formula corrects an oversight in \cite{W2}, p.~126 (the
first formula on p.~126 holds only if $x\ne 1$, so inversion hides one
term in the second).\vskip 12pt

The compatibility of the beta divisor with deletions (that is,
Proposition 3.9(2)) leads to the additivity property of the beta
invariant:
\proclaim{Corollary 3.12} If $e$ is not an isthmus, then $\beta(\Cal
L)=\beta(\Cal L-e)+\beta(\Cal L/e)$\endproclaim
\remark{Remark} Loops do not appear in this statement because our
matroids are loopless by assumption, cf.~the introduction.\endremark
\demo{Proof} If $e$ is not an isthmus, then $r(\Cal L)=r(\Cal L-e)$ so
$\alpha_*(\Til\omega_{\Cal L})=\Til\omega_{\Cal L-e}$ by Proposition
3.9(2). Using Propositions 3.3 and 3.10:
$$\align
(-1)^{r(\Cal L-e)}\beta(\Cal L-e) &=\gamma_0'\cdot\Til\omega_{\Cal
L-e}\\
&=(\gamma_0+\gamma_e)\cdot\Til\omega_{\Cal L}\qquad\text{by the
projection formula}\\
&=(-1)^{r(\Cal L)}\beta(\Cal L)+(-1)^{r(\Cal L/e)}\beta(\Cal
L/e)\quad.
\endalign$$
Since $r(\Cal L-e)=r(\Cal L)=r(\Cal L/e)+1$, the statement follows.
\qed\enddemo

If $e$ is an isthmus, an extra $-H_1'$ term appears in $\alpha_*(\Til
\omega_{\Cal L})$; if $\Cal L\ne [0,e]$, the argument in this proof
gives $\beta(\Cal L-e)=-\beta(\Cal L)+\beta(\Cal L/e)$ (since in this
case $r(\Cal L-e)=r(\Cal L/e)=r(\Cal L)-1$); and since $\Cal L-e=\Cal
L/e$ if $e$ is an isthmus, it follows that $\beta(\Cal L)=0$ in this
case. If $\Cal L=[0,e]$ itself is an isthmus, then $e=1$ and the extra
$H_1'$ term kicks in, giving $\beta([0,e])=1$ as it should
(cf.~\cite{W2}, 7.3.1(b)).

The vanishing of the beta invariant in the presence of an isthmus is a
particular case of the fact that the invariant vanishes on direct
sums. This will follow in a moment from Corollary 3.14 below; it could
also be checked easily by studying the deletion map $V_{\Cal
L_1\times\Cal L_2} @>>> V_{\Cal L_1}$. We leave this as a pleasant
exercise to the reader (although the conventional proof, which simply
takes the derivative of a product, is much easier!)\vskip 12pt

Proposition 3.9(3) translates into:
\proclaim{Corollary 3.13} If $x\in\Cal L$ is modular, and $y<x$, then
$$(-1)^{r(x)-r(y)}\beta([y,x])=\sum_{z\wedge x=y}B(z)$$
\endproclaim
\demo{Proof} Writing $\Til\omega$ modulo $F$ and using (3) from
Proposition 3.9 yields
$$i^*(\sum_{z\in\Cal L}B(z)\,H_z)=\sum_{0\le y\le x}B(y)''\,H_y''
-(r(1)-r(x))H_x''\quad.$$
But Proposition 3.5 says
$$\align
i^*(\sum_{z\in\Cal L}B(z)\, H_z)&=\sum_{z\in\Cal L}B(z)\, H_{z\wedge
x}''\\
&=\sum_{0\le y\le x}\left(\sum_{z\wedge x = y} B(z)\right)\,H_y''\quad;
\endalign$$
comparing the two expression and dotting with $\gamma_y''$ for $y<x$
gives the statement.\qed\enddemo

Setting $y=0$ in this corollary and isolating the term $B(0)$ gives:
$$(-1)^{r(1)} \beta(\Cal L)=(-1)^{r(x)}\beta([0,x])- \sum_{z\ne
0,z\wedge x=0}B(z)\tag*$$
if $x\ne 0$ is modular. The following statement follows:
\proclaim{Corollary 3.14} If $x\in\Cal L,x\ne 0$ is modular, then
$$\beta(\Cal L)=(-1)^{r(1)-r(x)}\beta([0,x])\sum_{y\wedge
x=0}\mu(0,y)$$
\endproclaim
\demo{Proof} Induction: if $r(\Cal L)=r(x)$ then $x=1$ and there is
nothing to prove; next, the terms in the summation in (*) are (up to
sign) beta invariants of lattices of lower rank, so we may apply the
statement to them (because $z\vee x$ is modular in $[z,1]$ and
$[z,z\vee x]\cong [0,x]$ if $z\wedge x=0$); doing so yields the
induction step.\qed\enddemo

The statement of the last corollary is a `modular decomposition'
expression for the beta invariant. It could also be derived easily
from Stanley's modular factorization theorem; the above proof,
however, seems more direct. For $x=e$ a rank@-1 element of $\Cal L$
(thus automatically modular), the corollary says
$$\beta(\Cal L)=(-1)^{r(\Cal L)-1}\sum_{y\ngeq e}\mu(0,y)\quad,$$
that is 7.3.1(d) in \cite{W2}. For $x=(1,0)$ in $\Cal L_1\times\Cal
L_2$, Corollary 3.14 implies the vanishing of the beta invariant on
direct sums: indeed in this case $y\wedge x=0\iff y\in\{0\}\times \Cal
L_2$, so $\sum_{y\wedge x=0}\mu(0,y)=0$ if $\Cal L_2\ne 0$.\vskip 12pt

One last observation: the last formula can also be written
$$B(0)-\sum_{y\ge e}\mu(0,y)=0\quad,$$
which translates back into $(\Til\omega_{\Cal L}+\sum_{y\ge
e}E_y)\cdot\gamma_0=0$
by Proposition 3.10 and Lemma~2.2, or in fact into
$$(\Til\omega_{\Cal L}+\sum_{y\ge e}E_y)\cdot\gamma_z=0\qquad\text{if
$z\ngeq e$}$$
(once more by compatibility with contractions). In other words,
writing
$$\Til\omega_{\Cal
L}+\sum_{y\ge e}E_y=\sum_{z\in\Cal L}{a_z}H_z$$
modulo $F$, we find $a_z=0$ necessarily for all $z\ngeq e$. This also
has a pretty geometric explanation. By induction and compatibility
with contractions, it is enough to show $\sum_{z\ngeq e}a_z=0$. Now
consider the hypersurface $D$ obtained by taking the union of a
general representative of $H_0$ and of all $E_y$ with $y\ge e$. $D$ is
non@-singular along $H_0$ and $E_e$ away from certain divisors of
these latter; the complements $\Til H_0$, $\Til E_e$ of these divisors
in $H_0$, $E_e$ are isomorphic to the complement of sets of
codimension at least 2 in $\P^{n-1}$: so their $\Pic$ is $\Z$,
generated by a hyperplane class $h$, and their canonical divisor is
$-n\,h$. Now
$$\omega_{\Cal L}+H_0+\sum_{y\ge e}E_y=\omega_{\Cal L}+D$$
restricts, by adjunction, to the canonical divisors of $\Til H_0$,
$\Til E_e$ on each of these: so restricting $\Til\omega_{\Cal L}+
\sum_{y\ge e}E_y = \omega_{\Cal L}+D+(n-r(1)) H_1$ to $\Til H_0,\Til
E_e$ and reading the coefficient of $h$ gives respectively
$$\gather
\sum_{z\ge 0}a_z = -n+(n-r(1))= -r(1)\quad,\\
\sum_{z\ge e}a_z = -n+(n-r(1))= -r(1)
\endgather$$
($H_0$ meets all $H_z$, while $E_e$ only meets the $H_z$ with $z\ge
e$). Comparing the two expression gives $\sum_{z\ngeq e}a_z=0$, as
needed.


\heading \S 4. Degrees of matroid varieties; Segre classes\endheading

The line bundles associated with the divisors
$$S(m)= m^{r(1)} H_0-\sum_{x\in\Cal L,x\ne 0} m^{r(1)-r(x)} E_x$$
introduced in \S 2 on the $n$-dimensional matroid variety $V^n=V_{\Cal
L}$ define for $m>0$ rational maps
$$\sigma_{m,n}:V^n\dashrightarrow\P^{N(m,n)}$$
to a projective space. We will write $\sigma_m$ for short
(disregarding $n$) because these maps are compatible with the natural
inclusions $V^n\subset V^{n+1}\subset\cdots$ discussed in \S 2.1
(since the $S(m)$ are).
\example{Example} For $m=1$ we have $S(1)= H_1$ modulo $F$,
and it follows that $\sigma_1$ is the blow@-up map $V^n @>>> \P^n$
followed by the projection with center $1\in\Cal L$.\endexample

Now define for $m>0,n>r(\Cal L)$:
$$d(m,n)=(\deg\,\sigma_m)\,(\deg\,\overline{\sigma_m(V^n)})$$
So $d(m,n)=0$ if $\dim \sigma_m(V^n)<n$, while $d(m,n)$ is just the
degree of $\overline{\sigma_m(V^n)}$ if $\sigma_m$ is generically
injective; for example $d(1,n)=0$ for all $n$.

At this stage we do not know a general formula for $d(m,n)$. In a
sense that is not surprising because, as we will show in a moment, the
characteristic polynomial of the original matroid can be recovered
from a fraction of the information carried by the $d(m,n)$'s. More
precisely, let $\{a\}_n$ denote the smallest nonnegative residue of
$a$ modulo $n$; then Theorem 4.4 will imply:

\noindent {\it Let $d(m,n)$ be the numbers defined above for the cycle
matroid of a simple graph $G$; and let $c$ be the number of components
of $G$. Then the value of the chromatic polynomial of $G$ at $m>0$
equals
$$m^c\,\{d(m,n)\}_n\quad,$$
where $n$ is an arbitrary sufficiently large prime.}\vskip 12pt

Also, observe that the $V^n$'s are birational to $\P^n$ (via the
blow@-up map), so that $\sigma_m$ and the $d(m,n)$ could be defined
starting from the original $\P^n$ in which $\Cal L$ is embedded, thus
bypassing the blow@-up construction. The right language to express
this is that of {\it Segre classes:\/} we will show that the $d(m,n)$
are determined by the Segre classes of specific subschemes of $\P^n$
supported on $1\in\Cal L$. Advances in the theory of Segre classes
could thus be relevant to problems of graph coloring!\vskip 12pt

In this section we say for short that a matroid is `nice' if it
belongs to the class defined in Proposition 2.5: that is, if all its
geometric contractions have characteristic polynomials with
nonnegative value at positive integers. In particular, for nice
matroids the divisor
$$\overline S(m)=\sum_{x\in\Cal L}p(\Cal L/x,m)\,H_x$$
introduced in \S 2 is generated by global sections for all $m>0$. Thus
graphical matroids, for example, are nice in this sense.
\proclaim{Lemma 4.1} For nice matroids\roster
\item the $\sigma_m$'s are in fact regular maps;
\item the pull@-back of the hyperplane class via $\sigma_m$ is
$\overline S(m)$.
\endroster\endproclaim
\demo{Proof} $\overline S(m)=S(m)$ modulo $F$: thus the rational maps
defined by $\overline S(m)$ and $S(m)$ agree on a non@-empty open
subset of $V^n$ (the complement of the $F$ divisors), hence they are
the same. Now $\overline S(m)$ is globally generated for nice matroids
and $m>0$, so the map is regular and $\overline S(m)$ is the
hyperplane section.\qed\enddemo

(2) implies:
\proclaim{Corollary 4.2} For nice matroids: $d(m,n)=\overline S(m)^n$
(the $n$-th self@-intersection of $\overline S(m)$ in
$V^n$).\endproclaim

Now computing $\overline S(m)^n$ is a challenge. The following trivial
observation is our only tool:
\proclaim{Lemma 4.3} $H_0^n=1$; $H_x^n=0$ for $x\ne 0$.\endproclaim
\demo{Proof} $H_0$ is the pull@-back of the hyperplane from $\P^n$ via
the blow@-up maps, so the first formula follows from the projection
formula. The second follows from the remark following Proposition 2.4:
the intersection of $n+1-r(x)\le n$ general representatives of $H_x$
is empty.\qed\enddemo

Still, this is enough to obtain the result mentioned in the
introduction:
\proclaim{Theorem 4.4} If $\Cal L$ is the lattice corresponding to a
nice matroid (e.g., a graphical matroid) and $n\ge r(\Cal L)$ is a
prime number, then
$$p(\Cal L, m)\equiv d(m,n)\pmod n$$
\endproclaim

In particular, let $\{a\}_n$ denote the smallest nonnegative residue
of $a$ modulo $n$. Then:
\proclaim{Corollary 4.5} If $\Cal L$ corresponds to a nice matroid,
$$p(\Cal L,m)=\{d(m,n)\}_n\qquad\text{for all primes $n\gg 0$}$$
\endproclaim

For graphs, the corollary implies the statement in italics given
earlier in this sections, by the relation between the chromatic
polynomial of a graph and the characteristic polynomial of its cycle
matroid.

\demo{Proof of the Theorem} From Corollary 4.2 and Theorem 2.3
$$\alignat 2
d(m,n)&=\overline S(m)^n=(\sum_{x\in\Cal L}p(\Cal L/x,m)\, H_x)^n &&\\
&\equiv \sum_{x\in\Cal L} (p(\Cal L/x,m)\, H_x)^n \pmod n&&\qquad\text
{since $n$ is prime}\\
&\equiv p(\Cal L,m)^n \,\,\pmod n &&\qquad\text{by Lemma 4.3}\\
&\equiv p(\Cal L,m) \quad\pmod n &&\qquad\text{by Fermat's little
theorem.}\qed
\endalignat$$
\enddemo

The $d(m,n)$ can alternatively be obtained in terms of the Segre
classes of subschemes of $\P^n$ supported on $1\in\Cal L$, whose
definition can be given without reference to the rest of the
construction. For each $m>0$, consider the subscheme $X(m,n)$ of
$\P^n$ defined by the intersection of all degree@-$m^{r(1)}$
hypersurfaces satisfying the multiplicity prescription mentioned in
the end of \S1---in short, multiplicity $m^{r(1)-r(x)}$ along $x$ for
all $x\ne 0$ in $\Cal L$ (such hypersurfaces do exist for nice
matroids: map the general representative of $\overline S(m)$ down to
$\P^n$; and their intersection is clearly supported on the maximal
subspace $1\in\Cal L$). Next, let $s_0(m,n)$ be the {\it degree of the
zero-dimensional component of the Segre class\/} $s(X(m,n),\P^n)$ (see
\cite{F}, Chapter 4, for the notion and properties of Segre classes).
The result is then

\proclaim{Theorem 4.6} If $\Cal L$ is the lattice corresponding to a
nice matroid (e.g., a graphical matroid) and $n\ge r(\Cal L)$ is a
prime number, then
$$p(\Cal L, m)\equiv m^{r(\Cal L)}-s_0(m,n)\pmod n$$
\endproclaim

Thus, the characteristic polynomial can be recovered in terms of these
numbers as well. Also, we note explicitly that the statement of the
four@-color@-theorem translates into:

{\it for a planar graph $G$, there exists a prime $n$ such that
$s_0(4,n)\not\equiv 4^r\pmod n$, where $r$ is the number of edges in a
spanning forest of $G$\/.}

Theorem 4.6 follows from the following relation between the $d(m,n)$
and the above Segre classes:
\proclaim{Lemma 4.7} For $m>0$ and $n>r(1)$:
$$d(m,n)=\left(m^{r(1)}\right)^n-\int_{X(m,n)}(1+m^{r(1)} H)^n\cap
s(X(m,n),\P^n)$$\endproclaim

\noindent where $H$ is the hyperplane class in $\P^n$, and $\int$
denotes degree in the sense of \cite{F},~\S1.4.
\demo{Proof} The linear system defined by the hypersurfaces of $\P^n$
satisfying the multiplicity prescription defines a rational map $\P^n
\dashrightarrow \P^{N(m,n)}$. Now we claim that this map, composed
with the blow@-up sequence defining $V^n$, gives the map $\sigma_m$
defined at the beginning of this section: this follows from Lemma 4.1,
since the proper transform of the hypersurfaces has class $\overline
S(m)$ in $V^n$. Then applying Proposition 4.4 in \cite{F} gives the
statement.\qed\enddemo

To prove Theorem 4.6, just read the Lemma modulo $n$ and apply Theorem
4.4.

Here is a table of $s_0(m,n)$ for the complete graph on three vertices:
\siltable
&&$s_0(m,n)$&& $m=2$ && 3 && 4 && 5 &\cr\tablerule
&& $n=3$  && 10 && 58 && 160 && 334 &\cr\tablerule
&&   4    && 30 && 213 && 726 && 1821 &\cr\tablerule
&&   5    && 74 && 692 && 3020 && 9308 &\cr\tablerule
&&   6    && 166 && 2143 && 12226 && 46795 &\cr\tablerule
&&   7    && 354 && 6510 && 49080 && 234282 &\cr.

We believe the $s_0(m,n)$ might in general be easier to control than
the $d(m,n)$.

To conclude, we mention that yet another congruence result similar to
Theorems 4.4, 4.6 can be stated in terms of Fulton's {\it canonical
classes\/} (4.2.6(a) in \cite {F}). For this, denote by $X_H(m,n)$ the
general hyperplane section of $X(m,n)$, and by $c_0(m,n)$ the degree
of $c_0(X(m,n))-c_0(X_H(m,n))$ (notations as in \cite{F}, Example
4.2.6); then one can show
$$c_0(m,n) \equiv s_0(m,n)\pmod n$$
for $n$ prime. Unfortunately, few properties of Fulton's canonical
classes are known as yet.


\enddocument